\documentclass[10pt]{amsart}

\usepackage{amsmath}
\usepackage{amssymb}
\usepackage{amscd}
\usepackage{subfigure}
\usepackage{graphicx}

\newcommand{\N}{\mathbb N}

\newcommand{\R}{{\mathbb R}}
\newcommand{\Q}{{\mathbb Q}}
\newcommand{\Z}{{\mathbb Z}}
\newcommand{\C}{{\mathbb C}}
\newcommand{\Ak}{{\mathcal A}}

\newtheorem{theorem}{Theorem}[section]
\newtheorem{lemma}[theorem]{Lemma}
\newtheorem{prop}[theorem]{Proposition}
\newtheorem{cor}[theorem]{Corollary}

\numberwithin{equation}{section}

\begin{document}
\title{A family of non-sofic beta expansions}
\address{Institute of Mathematics, University of Tsukuba\\ 
Tennodai 1-1-1, Tsukuba, Ibaraki\\
305-8571 Japan}
\email{akiyama@math.tsukuba.ac.jp}
\author{Shigeki Akiyama}
\date{}
\thanks{
This research was supported by the
Japanese Society for the Promotion of Science (JSPS), Grant in aid
21540012.}
\maketitle

\begin{abstract}
Let $\beta_n>1$ be a root of $x^n-x-1$ for $n=4,5,\dots$. We will prove that
$\beta_n$ is not a Parry number, i.e., the associated beta transformation 
does not correspond a sofic symbolic system. A generalization is shown in the last 
section.
\end{abstract}


\section{Beta expansions}
Fix a real number $\beta>1$. The map from $[0,1)$ to itself defined by
$T_{\beta}(x)=\beta x-\lfloor \beta x \rfloor$ 
is called the {\it beta transformation}. Putting $a_n=\lfloor \beta T^{n-1}_{\beta} (x)\rfloor$, we obtain an expansion:
$$
x = \frac {a_1}{\beta}+ \frac {a_2}{\beta^2}+\dots
$$
with $a_i\in \Ak:=\Z \cap [0,\beta)$, which gives a generalization of decimal expansion 
to the real base $\beta$. Let $\Ak^{\N}$ (resp. $\Ak^{\Z}$) 
be the set of right infinite (resp. bi-infinite) words over $\Ak$ 
which is compact by the product topology of $\Ak$. 
Define $d_{\beta}:[0,1)\rightarrow \Ak^{\N}$ by $d_{\beta}(x)=a_1a_2\dots$. 
The {\it expansion of one} of $\beta$
is the infinite word $c_1c_2\dots\in \Ak^{\N}$ obtained 
as a limit of the expansion $1-\epsilon$ when $\epsilon>0$ tends to zero, 
which is denoted by $d_{\beta}(1-0)$. The map $d_\beta$ is not surjective and 
the image $d_{\beta}([0,1))$ is characterized as
$$
\{ \xi=(\xi_n)\in \Ak^{\N} \ |\ s^n(\xi)\ll d_{\beta}(1-0)\quad (n=0,1,\dots) \}
$$
where $s$ is a shift operator $s((\xi_n))=(\xi_{n+1})$, and 
$\ll$ is the natural lexicographic order on $\Ak^{\N}$. We say that
$\xi\in \Ak^{\N}$ is {\it admissible} if it satisfies the {\it Parry condition}
$$
s^n(\xi)\ll d_{\beta}(1-0)\quad (n=0,1,\dots),
$$
see \cite{Parry:60, Ito-Takahashi:74}.
Let $\Ak^*$ be the set of finite words over $\Ak$. An element $w\in \Ak^*$
is admissible if $w0^{\infty}=w00\dots$ is admissible.  Define 
a compact subset of $\Ak^{\Z}$ by
$$
X_{\beta}=\{ (\xi_n)\in \Ak^{\Z}\ |\ \xi_n\xi_{n+1}\dots \xi_m \text{ is admissible for all } n \text{ and } m \text{ with } n<m \}.
$$
The symbolic dynamical system $(X_{\beta},s)$ is called {\it beta shift}.
We see that $(X_{\beta},s)$ is a subshift of finite type if and only if
$d_{\beta}(1-0)$ is purely periodic. Further $(X_{\beta},s)$ is sofic
if and only if $d_{\beta}(1-0)$ is eventually periodic. 
We say that $\beta$ is a
{\it simple Parry number} if $(X_{\beta},s)$ is a shift of finite type, and a
{\it Parry number}\footnote{Parry coined it {\it beta number} but it is  
confusing to say $\beta$ is a beta number. Recent articles use this name.} 
if $(X_{\beta},s)$ is sofic. It is well known that
$\beta$ is sofic if $\beta$ is a Pisot number, that is, a real algebraic integer 
greater than one whose all conjugates lie within the open unit disk.
In fact, this follows from a
general fact that beta expansions of elements of $\Q(\beta)\cap [0,1)$ are
eventually periodic provided $\beta$ is a Pisot number \cite{Bertrand:77,Schmidt:80}. 
Boyd \cite{Boyd:89, Boyd:96_2}
showed that Salem numbers of degree 4 are Parry numbers, and gave some 
heuristic discussion on the existence of non-Parry Salem number of higher degree. 
However until now, 
we have no idea to prove that $d_{\beta}(1-0)$ is not eventually periodic
when $\beta$ is a Salem number.
In this note, we will show the following 

\begin{theorem}
\label{Main}
Let $\beta_n>1$ be the root of $x^n-x-1$ for $n=2,3,\dots$. Then
$\beta_n$ is a Parry number if and only if $n=2,3$.
\end{theorem}

According to \cite{Parry:60}, we know that if $\beta$ is a Parry number, then it must be
a Perron number whose conjugates has modulus less than 2. Here a Perron number is
an algebraic integer greater than one whose all other 
conjugates has modulus strictly less than the number itself. 
Solomyak \cite{Solomyak:94}
further 
studied
distribution of conjugates of Parry numbers, describing the intriguing 
region $\Phi$ where the conjugates densely lie. This 
improves the modulus bound to $(1+\sqrt{5})/2$. He also gave an example
of a non-Parry Perron number $(1+\sqrt{13})/2$
whose conjugate lie in the interior of $\Phi$.
Theorem \ref{Main} seems to be the first result on a family of 
non-Parry Perron numbers whose conjugates lie\footnote{
Seemingly the conjugates of $\beta_n$ 
are in the interior of $\Phi$ from Fig. 2 in \cite{Solomyak:94}.
However it is not so easy to show this, because of the fractal nature
of the boundary of $\Phi$.} in $\Phi$.

The key of the proof is {\it Lagrange inversion formula} which gives 
the inverse 
of Taylor expansion of a holomorphic function defined in some region.
As Theorem \ref{Main} covers all $n$, 
we must rely on numerical computation. 
The dependencies to computer are sketched along the proofs. 
If we permit a finite number of exceptions, then 
the proof becomes computer independent
and we can treat wider cases. 
A generalization of Theorem \ref{Main} in this sense
is given in the last section.

Hereafter Landau $O$ symbol will be in abusive use :
$f(x)=O(g(x))$ means there exists a constant $C$ that 
$|f(x)|\le C|g(x)|$ for all $x$ in an appropriate ball (possibly centered at $\infty$)
which is clear from the context. Vinogradoff symbols are 
not used. We write $n \gg 1$ only to mean that $n$ is sufficiently large. 

\section{Proof}

Let $\beta$ be a non-Pisot Perron number. Then one can select 
a conjugate $\gamma$ of $\beta$ that $|\gamma|>1$.
Let $x'$ be the image of $x$ by the conjugate map from
$\Q(\beta)$ to $\Q(\gamma)$ and $d_{\beta}(1-0)=c_0c_1\dots$.
Put
$$
T_{\beta}^{k}(1-0)=\beta^k\left(1-\sum_{m=1}^k \frac {c_m}{\beta^m}\right) 
\in \Z[\beta].
$$
Note that $T_{\beta}^{0}(1-0)=1$ and we have
$$
T_{\beta}^{k}(1-0)=\sum_{m=1}^{\infty} \frac {c_{m+k}}{\beta^m}.
$$

\begin{lemma}
\label{Non}
If there is $k\in \N$ with $|(T_{\beta}^k(1-0))'|> \frac {\lfloor \beta\rfloor}{|\gamma|-1}$,
then $\beta$ is not a Parry number.
\end{lemma}

\proof
Putting $x_m=T_{\beta}^m(1-0)$, we have $x_{m+1}=\beta x_m-c_{m+1}$.
Since $|x_k|> \lfloor \beta \rfloor/(|\gamma|-1)$, we have
$$
|x_{m+1}'|= |\gamma x_m' -c_{m+1}|> |x_m'|
$$
for $m\ge k$. Therefore the sequence $(|x_m'|)_{m=1,2,\dots}$ diverges, which is
impossible if $c_i$ is eventually periodic. 
\qed
\bigskip

This lemma gives a computational way to show that $\beta_n$ is not 
a Parry number for a fixed $n$.

For $n=2$, $\beta_2=(1+\sqrt{5})/2$ is the most known Pisot number, the golden mean. 
It is also well known that $\beta_3$ the smallest Pisot number \cite{Smyth:71,
Pisot-Salem:92}. 
We will show that $\beta_n$ for $n\ge 4$ is not a Parry number. 

\begin{lemma}
\label{Perron}
$\beta_n\ (n\ge 4)$ is a non-Pisot Perron number.
\end{lemma}

\proof
Let $V=\{1,\dots, n\}$ and define the directed edge $E$ by
$$
i\rightarrow i+1 \quad (i=1,2,\dots,n-1), n\rightarrow 1, n\rightarrow 2.
$$
The adjacency matrix of this graph is clearly primitive and its
Perron-Frobenius root is $\beta_n$, which shows that 
$\beta_n$ is a Perron number.
From $(\beta_n)^{n+1}-\beta_{n}-1=\beta_n^2-1>0$, we see
$$
\beta_2>\beta_3>\beta_4> \dots >1
$$
Since $\beta_3$ is the smallest Pisot number, $\beta_n$ for $n\ge 4$ is
not a Pisot number.
\qed
\bigskip

\begin{lemma}
\label{Irred}
The polynomial $x^n-x-1$ is irreducible over $\Q$ for $n\ge 2$.
\end{lemma}

\proof
This result is due to Selmer \cite{Selmer:56}. 
\qed
\bigskip

B\"urmann-Lagrange formula 
is discussed in Part~I-Chap.~7 of \cite{Hurwitz-Courant:44}.
We briefly review it in a special form, to obtain 
an explicit truncation error bound.
Denote by $B(x,r)$ the ball of radius $r$ centered at $x$.
Let $g(z)$ be a holomorphic function
with $g(0)=0$ and $g'(z)\neq 0$ in $z\in B(0,r)$.
Then $g$ is locally univalent and admits a holomorphic inverse 
which is to be made explicit.
Define a function
$$
h(w)=\frac {1}{2\pi \sqrt{-1}} \oint_C \frac{\zeta g'(\zeta)}{g(\zeta)-w} d\zeta
$$
where $C$ is the counter-clockwise contour which circumscribes $B(0,r)$.
Since $g'(z)$ does not vanish, 
by residue theorem we have $h(g(z))=z$ 
in a neighborhood of the origin, and hence in $B(0,r)$ by identity theorem on
holomorphic functions. 
Using 
$$
\frac 1{1-z}= \sum_{k=0}^m z^k + \frac {z^{m+1}}{1-z}
$$
we have
\begin{equation}
\label{Inv}
h(w)= c_1w+\dots + c_m w^m + 
\frac {1}{2\pi \sqrt{-1}} \oint_C \frac{\zeta g'(\zeta) w^{m+1}}{g(\zeta)^{m+1}
(g(\zeta)-w)} d\zeta
\end{equation}
with
$$
c_k=\frac {1}{2\pi n \sqrt{-1}} \oint_C \frac{d\zeta}{g(\zeta)^{k}}.
$$
This (\ref{Inv}) is the 
{\it Lagrange inversion formula} in a complex analytic form. 
A different formulation is found in pp.131-133 of
\cite{Whittaker-Watson:27}. It has many interesting applications in combinatorics. 

\begin{prop}
\label{Roots}
Fix $m\in \Z$. For an integer $n\ge 12|m|$, 
there is a root of $x^n-x-1$ which satisfies the asymptotic formula:
$$
\exp\left(\frac {2\pi m\sqrt{-1}}n\right) +\frac{\log 2}{n}+\frac{
(1+\log 2)\log 2 + 2 \sqrt{-1}\pi m(1+ \log 4)}{2 n^2}+ C(n)
$$
with $C(n)=O\left(\frac{1}{n^3}\right)$. 
\end{prop}

\begin{figure}[h]
\includegraphics[height=6cm]{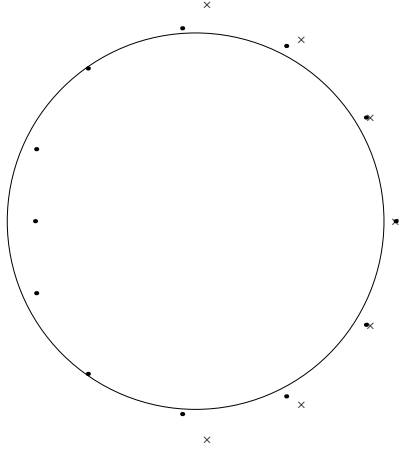}
\caption{Roots of $x^{12}-x-1$ (black dots) 
and approximations ($\times$ dots) by the formula of Proposition \ref{Roots}. 
\label{Approx}}
\end{figure}

\proof
Consider a root $\gamma$ of $x^n-x-1$ lying in a ball $B(1,1/2)$. 
Since $|\arg \gamma|<\pi/6$, we have
$$
\frac 1n=\frac {\log(\gamma)}{\log (1+\gamma)+2\pi m \sqrt{-1}}
$$
in the principal branch of logarithm,  $m\in \Z$ and $|m|<n/12$.
We fix $m$ and study the asymptotic behavior of $\gamma$ when $n$ tends to $\infty$.
Introduce a complex variable $z=\gamma-1$ to define
$$
g(z)=\frac {\log(z+1)}{\log (z+2)+2\pi m \sqrt{-1}}.
$$
Then $g(z)$ is holomorphic $g(0)=0$ and $g'(z)\neq 0$ in $B(0,1/2)$.
Lagrange inversion (\ref{Inv}) gives
$$
h(w)= (\log 2+ 2\pi m \sqrt{-1}) w + \left(\frac{(1+\log 2)\log 2}2 + 
\sqrt{-1}\pi m(1+ \log 4)-2 \pi^2 m^2 \right) w^2 + E(w)
$$
with
$$
E(w)= \frac {1}{2\pi \sqrt{-1}} \oint_C \frac{\zeta g'(\zeta)w^3}{g(\zeta)^{3}(g(\zeta)-w)} d\zeta =O(w^3)
$$
where $C$ is the contour for $B(0,1/2)$. 

Putting $w=1/n$, the Taylor expansion of 
$\exp(2\pi m\sqrt{-1}/n)$ leads to the required asymptotic formula.
\qed

We see that $x^n-x-1$ has a unique root greater than $1$. Denote this root by
$\beta_n$. Let $\gamma_n$ be the complex root of $x^n-x-1$ 
closest to $\beta$ in $\C$ with $\Im \gamma_n>0$. 

\begin{cor}
\label{Est}
\begin{eqnarray}
\label{beta}
&&\left|\beta_n-\left(1+\frac{\log 2}{n}\right)\right| \le \frac{2}{3n^2} 
\qquad (n\ge 8)
\\
\label{gamma}
&&\left|\gamma_n-\left(1 +\frac{\log 2+2\pi\sqrt{-1}}{n}
\right)\right| \le \frac{24}{n^2} \qquad (n\ge 6).
\end{eqnarray}
\end{cor}
\bigskip

Note that Selmer \cite{Selmer:56} obtained a weaker form of (\ref{beta}).
In the course of the later 
proofs, we shall use numerical values of $\beta_n$ and $\gamma_n$
for small $n$'s. However they are not literally small. 
In particular, we will use
$\gamma_n$ with $n\le 3605$ which is computed by the complex Newton method
with the initial value $1 +(\log 2+2\pi\sqrt{-1})/n$.

\proof
We use $g, E_n$ in the proof of Proposition \ref{Roots}.  
For $m=0$, we use the numerical estimates
$\min\{|g(\zeta)| \ |\ |\zeta|=1/2\}\ge 0.44$ and 
$\max\{|g'(\zeta)|\ |\ |\zeta|=1/2\}\le 8$.
Assuming $n\ge 100$, it suffices to have
$$
\frac{(1+\log 2)\log 2}{2n^2} + \frac {8 \cdot 0.5^2}{0.44^3 n^3 \cdot (0.44-1/100)} 
< \frac 2{3n^2}.
$$
This is valid for $n\ge 684$. We can check the statement for $6\le n \le 683$ 
by numerical computation.
For $m=1$, we use 
$\min\{|g(\zeta)| \ |\ |\zeta|=1/2\}\ge 0.0636$ and 
$\max\{|g'(\zeta)|\ |\ |\zeta|=1/2\}\le 0.32$. Then the similar inequality 
$$
\frac{|(1+\log 2)\log 2+2(1+\log 4) \sqrt{-1} \pi-4\pi^2|}{2n^2}+
\frac {0.32 \cdot 0.5^2}{0.0636^3 n^3 \cdot (0.0636-1/1400)} 
< \frac {24}{n^2}.
$$
holds for $n\ge 1441$. Remaining $8\le n\le 1440$ are confirmed by 
direct computation.
\qed
\bigskip

We derive 
 three lemmas \ref{NextOne}, \ref{Ineq} and \ref{Abs}
 which are used in the proof of Theorem \ref{Main}.
Similarly to the proof of Corollary \ref{Est}, 
their proofs are finished for large $n$'s by (\ref{beta}) and (\ref{gamma}), while  
the remaining small $n$'s have to be checked by numerical computation. 

Since $\beta_n<2$ for all $n\ge 2$, we have $\Ak=\{0,1\}$ and
$c_1=1$. Let $m_0\ge 2$ the smallest index that $c_{m_0}=1$. 
First we have

\begin{lemma}
\label{NextOne}
$$
m_0\ge \frac{n \log n}{\log 2}
$$
for $n\ge 8$.
\end{lemma}

\proof
By the definition of $d_{\beta}(1-0)$, we have
$m_0= \left\lfloor \frac {\log(1-1/\beta_n)}{\log (1/\beta_n)} \right\rfloor$.
By (\ref{beta}), it suffices to show
$$
- \frac {\log \left(\frac {\log 2}n+ \frac 2{3n^2}\right)}{\log 
\left(1+ \frac {\log 2}n + \frac 2{3n^2}\right)} 
> \frac {n\log n}{\log 2}
$$
for $n\ge 8$.
\qed
\bigskip

More precise computation gives
$$
m_0=\frac{n \log n-n \log \log 2}{\log 2}
-\frac{\log n}{2 \log 2} + O(1),
$$
but we do not need this precision for the later use.

\begin{lemma}
\label{Ineq}
For $n\ge 6$ and $m_1\ge \frac {n\log n}{\log 2}$,  we have
$$
\left|\gamma_n^{m_1}(1-1/\gamma_n)\right|>4
$$
and 
$$
\left|\gamma_n^{m_1-2}\right|>\frac {n}2.
$$
\end{lemma}

\proof
Let $C$ be the counter-clockwise path around $B(0,1/2)$,
Taylor expansion 
$$
\log(1+z)=\sum_{i=1}^m \frac {(-1)^{i-1}z^m}i + \frac 1{2\pi \sqrt{-1}}
\oint_C \frac {\log(1+\zeta)z^{m+1}}{\zeta^{m+1}(\zeta-z)} d\zeta
$$
gives an estimate
$$
|\log(1+z)-z|\le  \frac{2\log 2}{1/2-|z|} |z^{2}|
$$
for $|z|<1/2$.
Since $|\gamma|>1$, we have
$$
\left|\gamma_n^{m_1}(1-1/\gamma_n)\right|
\ge \left|\gamma_n^{n\log n/\log 2-1}(\gamma_n-1)\right|.
$$
As
$$
\log(\gamma_n)=\log\left(1+\frac{\log 2+ 2\pi \sqrt{-1}}n+ \frac A{n^2}\right)
$$
for $|A|\le 24$, we have
$$
\log(\gamma_n)= \frac{\log 2+ 2\pi \sqrt{-1}}n+ \frac A{n^2} +
\frac {B}{n^2}
$$
with $|B|\le \frac {2\log 2}{1/2-7/2000}\cdot 6.4^2\le 115$ for $n\ge 2000$.
Here we used an estimate
$$
\left| \frac{\log 2+ 2\pi \sqrt{-1}}n + \frac A{n^2}\right|\le \frac {6.4}{n}
$$
valid for $n\ge 305$.
Therefore we have
\begin{equation}
\label{GamEst}
\log(\gamma_n)= \frac{\log 2+ 2\pi \sqrt{-1}}n+ \frac C{n^2}
\end{equation}
with $|C|\le 139$. Consequently
\begin{eqnarray*}
&&\left(\frac {n\log n}{\log 2} -1\right) \log (\gamma_n)\\
&=& \log n + \frac {2\pi \sqrt{-1} \log n}{\log 2}
+ \frac {C\log n}{n \log 2} - \frac {\log 2+2\pi \sqrt{-1}}n - \frac C{n^2}\\
&=& \log n + \frac {2\pi \sqrt{-1} \log n}{\log 2}
+ \frac {D\log n}{n}
\end{eqnarray*}
with $|D|\le 201$.
On the other hand, we have
\begin{eqnarray*}
\log(\gamma_n-1)&=&\log\left(\frac {\log 2+2\pi\sqrt{-1}}n +\frac A{n^2}\right)\\
&=& \log(\log 2 +2\pi\sqrt{-1})-\log n + \log\left(1+\frac A{n(\log 2+ 2\pi \sqrt{-1})}
\right)\\
&=& \log(\log 2+2\pi\sqrt{-1})-\log n + \frac A{n(\log2+ 2\pi \sqrt{-1})}+ \frac E{n^2}
\end{eqnarray*}
where
$$
|E|\le \frac {2\cdot 3.8^2 \log 2}{1/2-3.8/2000}\le 41.
$$
Here we used $|A/(\log 2+2\pi \sqrt{-1})|\le 3.8$.
Summing up, we have
$$
|\gamma_n^{n\log n/\log 2-1} (\gamma_n-1)|
= |\log 2 +2\pi \sqrt{-1}| \exp\left(\frac {D\log n}n +\frac Fn+\frac {E}{n^2}\right)
$$
with $|D|<201$, $|E|<41$, $|F|\le 3.8$ and $n\ge 2000$. 
For $n\ge 3606$, the last value exceeds $4$
and we obtain the first estimate of Lemma \ref{Ineq}. 
For $6\le n<3605$, we have to rely on the numerical computation.
For the second estimate, using (\ref{GamEst}),
\begin{eqnarray*}
&&\Re \left(\left(\frac {n\log n}{\log 2}-2\right) \log (\gamma_n)\right) \\
&&= \Re(\left(\frac {n\log n}{\log 2}-2\right)\left( 
\frac{\log 2+ 2\pi \sqrt{-1}}n+ \frac C{n^2}\right)\\
&&= \log n + \Re(C) \left(\frac {\log n}{n \log 2} - \frac 2{n^2} \right)
-\frac {2 \log 2}n \\
&&= \log n + G \frac {\log n}{n} 
\end{eqnarray*}
with $|G|\le 201$ and $n\ge 2000$. So we have
$$
|\gamma_n^{m_1-2}|\ge n \exp\left( G \frac {\log n}n \right) > \frac n2
$$
for $n\ge 2237$. The remaining $6\le n<2236$ is confirmed by numerical computation.
\qed

\begin{lemma}
\label{Abs}
For $n\ge 8$, we have
$$
\frac 1{\left|\gamma_n\right|-1} \le \frac {3n}2
$$
\end{lemma}

\proof
Using $(\ref{gamma})$, we have
$$
|\gamma \overline{\gamma}|
=1+ \frac {2\log 2}n+\frac {2\Re A}{n^2}+ \frac {|\log 2+2\pi \sqrt{-1}|^2}{n^2}
= 1+ \frac {2\log 2}n + \frac H{n^2}
$$
with $|H|\le 90$. We see
$$
\left| \sqrt{1+z}- \left(1+\frac z2\right)\right| \le \frac {\sqrt{6} |z|^2}{1/2-|z|},
$$
in a similar manner. Thus we obtain
\begin{equation}
\label{Recip}
|\gamma_n|-1 = \frac {\log 2}n + \frac H{2n^2} + \frac J{n^2}
\end{equation}
with $|J|\le \frac {1.5^2\sqrt{6}}{1/2-1.5/2000}\le 12$ for $n\ge 2000$.
Here we used an estimate
$$
\frac{2\log 2}n+\frac H{n^2}\le \frac{1.5}{n}
$$
for $n\ge 800$. Using (\ref{Recip}), we see that the statement is true for $n>2153$. 
The remaining $8\le n\le 2152$ are checked by direct computation.
\qed
\bigskip

{\it Proof of the Theorem \ref{Main}.}

Since every finite subword of $d_{\beta}(1-0)$ is admissible,
by the Parry condition, $10^{t}1\in \Ak^*$ is not admissible for $t<m_0-2$.
From the definition of $m_0$, we have
$c_{m_0+i}=0$ for $1\le i\le m_0-2$. By Lemma \ref{Non}, our goal is to prove
\begin{equation}
\label{Goal}
\left|(T_{\beta}^{2m_0-2}(1-0))'\right| > \frac 1{\left|\gamma_n\right|-1}.
\end{equation}
From Lemma \ref{Ineq} and
 $T_{\beta_n}^{2m_0-2}(1-0)=\beta_n^{2m_0-2}(1-\beta_n^{-1}-\beta_n^{-m_0})$,
 we have
\begin{eqnarray*}
\left|(T_{\beta}^{2m_0-2}(1-0))'\right|&=&\left|\gamma_n^{2m_0-2}(1-\gamma_n^{-1}-\gamma_n^{-m_0})\right|\\
&\ge& \left|\gamma_n^{2m_0-2}(1-\gamma_n^{-1})\right| - \left|\gamma_n^{m_0-2}\right|\\
&\ge& 3\left|\gamma_n^{m_0-2}\right| > 3n/2.
\end{eqnarray*}
which proves the theorem for $n\ge 8$
with the help of Lemma \ref{Abs}. For $n=6,7$, we can check (\ref{Goal})
directly. For $n=4$, we have 
$$
d_{\beta}(1-0)=10000 00010 00000 00000 01000 00000 10000 0\dots
$$
and 
$$
\left|(T_{\beta}^{m}(1-0))'\right| > \frac 1{\left|\gamma_n\right|-1}
$$
for $m=35$. For $n=5$, we get
$$
d_{\beta}(1-0)=10000 00000 00100 00000 00000 00\dots
$$
and one can take $m=26$.\qed

\section{A generalization}

There may be several ways to generalize Theorem \ref{Main}. 
Here we present a straight forward one. 

\begin{theorem}
\label{Gen}
Let $G$ be a polynomial with non negative integer coefficients
such that $G(1)>1$, $G(0)\neq 0$ and it is not a power of another polynomial.
Let $\alpha_n>1$ be the real root of
$
x^n - G(x).
$
Then there is a positive integer $n_0$ 
that $\alpha_n$ is a non-Parry Perron number for $n\ge n_0$.
\end{theorem}

\proof
Put $F(x)=x^n-G(x)$. Since $x>1$ implies $F'(x)>0$ for $n\gg 1$, 
$F(1)<0$ shows that there is a unique root $\alpha_n>1$ of $F$. Fixing $r>1$, 
from the non negativity
of the coefficients of $G$, we see that
$G(r)$ is the maximum of $|G(r\zeta)|$ for all $\zeta$ with
$|\zeta|=1$. It is unique in the sense that $|G(r\zeta)|=G(r)$ implies $\zeta=1$.
We know that $\alpha_n$ is a Perron number by virtue of 
Rouch\'e's theorem for a counter-clockwise circular path of 
radius $\alpha_n$ centered at $0$
avoiding outward the real root $\alpha_n$ by small perturbation.
Let $K(F)$ be the factor of $F$ whose leading coefficient 
is equal to the one of $F$, having 
properties that every root of $K(F)$ is
not a root of unity and $F/K(F)$ is a product of cyclotomic polynomials.
Theorem 5 of Schinzel \cite{Schinzel:65} reads that
there exists a positive integer $n_1$ that $K(F)$ is irreducible for
$n\gg 1$ and $(n,n_1)=1$.
Reviewing its proof, $n_1$ must be greater than one only when $x^n-G(y)$ is reducible
as a polynomial of $\Q(y)[x]$, which happens 
when $G(y)=h(y)^k$ with $k\ge 2$ or $G(y)=-4 h(y)^{4}$ for some $h\in \Q(y)$ 
by the theorem of Capelli (Theorem 9.1 in \cite{Lang:02}). Thus under
our assumption, we can take $n_1=1$.
The remainder of the proof proceeds similar to Theorem \ref{Main}.
Applying Lagrange inversion formula to
$$
g(z)=\frac {\log(z+1)}{\log G(z+1)+2\pi m \sqrt{-1}},
$$
we obtain the asymptotic expansion
$$
\alpha_n=1+\frac{\log G(1)}{n} + O\left(\frac 1{n^2}\right)
$$
and find a conjugate 
$$
\eta_n=1 +\frac{\log G(1)+2\pi m \sqrt{-1}}{n} + O\left(\frac 1{n^2}\right)
$$
for $n\gg 1$.
We select $m\in \N$ with $\exp(2\pi m/\sqrt{3})>G(1)$.
Clearly 
$\alpha_n$ and $\eta_n$ are the roots of $K(F)$
for $n\gg 1$. We obtain
asymptotic expansions:
$$
m_0:= \left\lfloor \frac {\log (1-1/\alpha_n)}{\log (1/\alpha_n)} \right\rfloor
= \frac {n\log n}{\log G(1)} -\frac {n \log \log G(1)}{\log G(1)} + O(\log n),
$$
$$
|\eta_n^{2m_0-2}(1-1/\eta_n)|= \frac {|\log G(1)+2\pi m \sqrt{-1}|}{(\log G(1))^2} n
+ O(\log n),
$$
$$
|\eta_n^{m_0-2}|= \frac n{\log G(1)} + O(\log n)
$$
and
$$
\frac 1{|\eta_n|-1} = \frac n{\log G(1)}+ O(1).
$$
Therefore
\begin{eqnarray*}
|T^{2m_0-2}(1-0)'| &\ge& |\eta_n^{2m_0-2}(1-1/\eta_n)|- |\eta_n^{m_0-2}|\\
&=& \frac {|\log G(1)+2\pi m \sqrt{-1}|}{(\log G(1))^2}n - \frac{n}{\log G(1)} + O(\log n)\\&>& \frac 1{|\eta_n|-1} 
\end{eqnarray*}
The last inequality holds for $n\gg 1$ by the choice of $m$.
\qed
\bigskip

We may expect some generalization
of Theorem \ref{Gen} for polynomials of the form
$x^n f(x)-g(x)$ for fixed $f$ and $g$, as
Lagrange inversion formula likewise applies.

\begin{figure}[h]
\begin{center}
\subfigure[$x^3-x^2+2x+2$]{
\includegraphics[height=3.5cm]{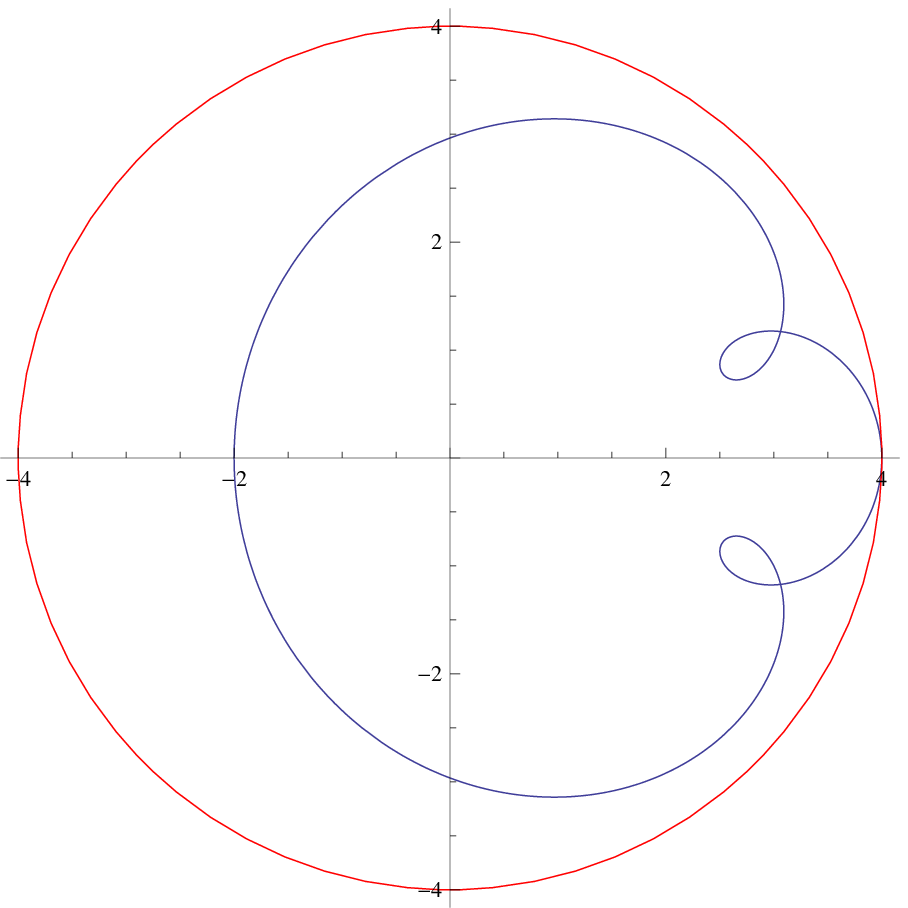}
}
\hspace{0.5cm}
\subfigure[$x^3+3x^2-x+1$]{
\includegraphics[height=3.5cm]{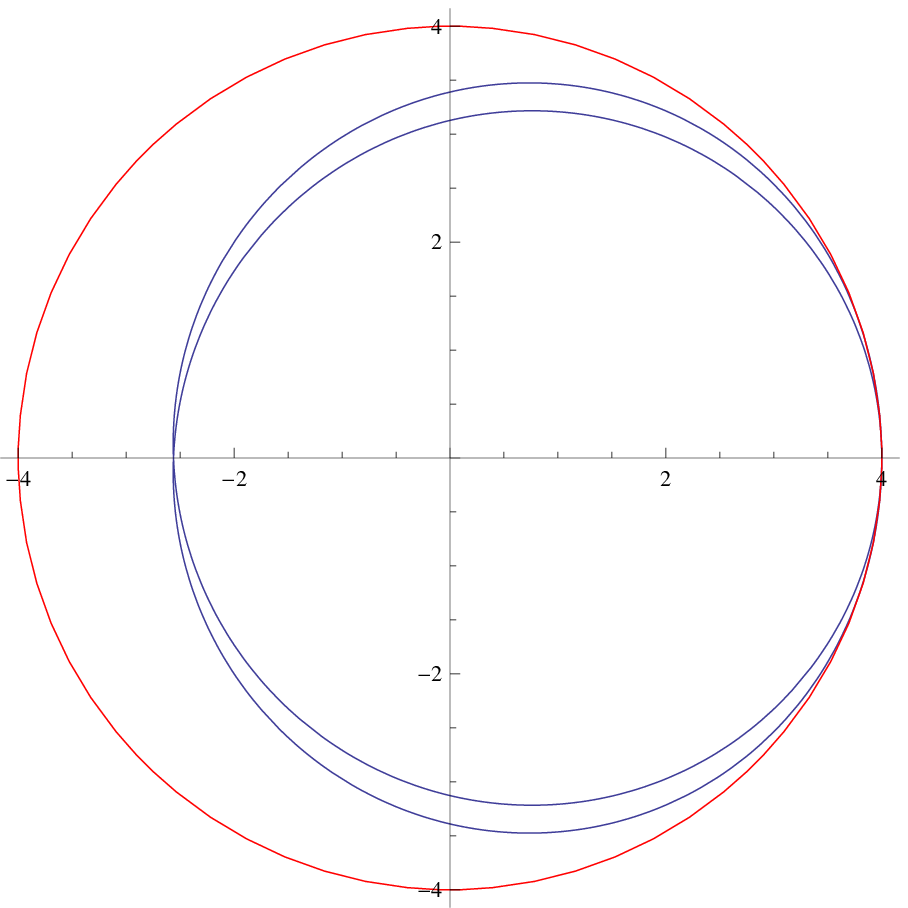}
}
\hspace{0.5cm}
\subfigure[$-x^3+3x^2+x+2$]{
\includegraphics[height=3.5cm]{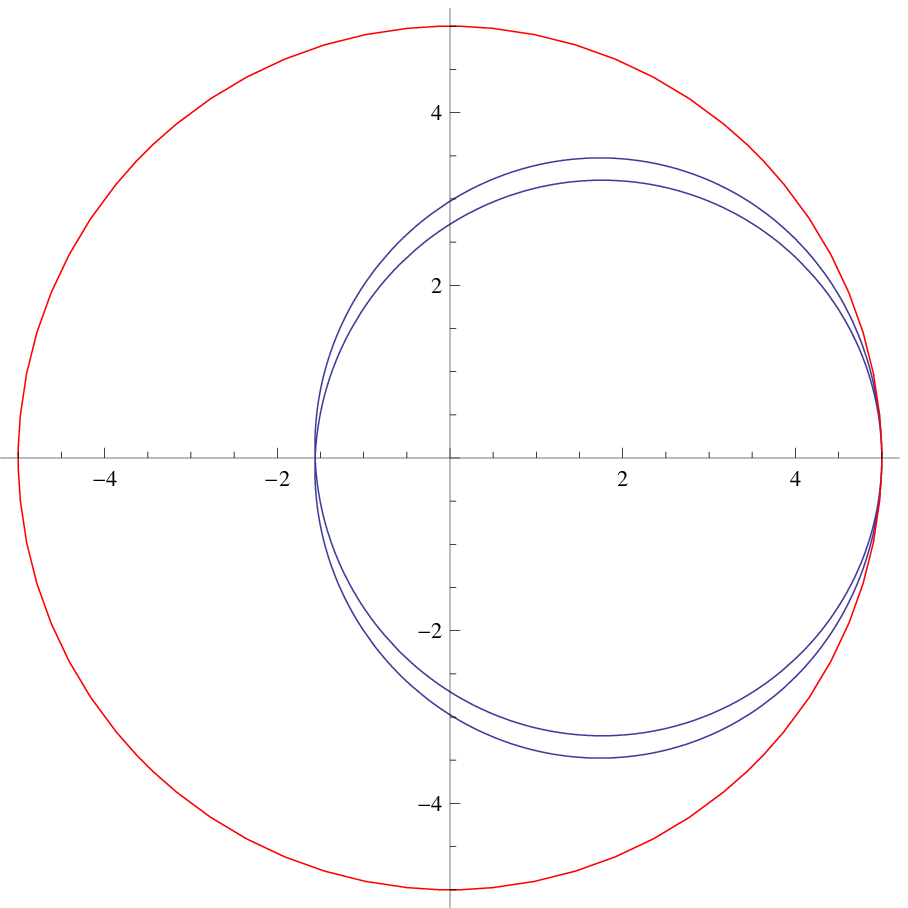}
}
\end{center}
\caption{Curves for $G(\exp(\sqrt{-1} t))$ and a circle of radius $G(1)$
\label{PerronCheck}}
\end{figure}

Without any change of the proof, 
the non-negativity condition of coefficients of $G$ can be relaxed to:
$$
\exists r_0>1,\ 1< \forall r<r_0,\ \forall \zeta\neq 1 \text{ with } |\zeta|=1 
\qquad |G(r\zeta)|< G(r).
$$
This is a geometric condition on a surface $G(r\exp(t\sqrt{-1}))$
parametrized by $r$ and $t$, which seems hard to check, 
but fulfilled by $G(x)=x^3-x^2+2x+2$, for e.g.
This is confirmed by checking the condition in the limit case $r=1$ 
(see Figure~\ref{PerronCheck}), 
and the fact that the surface is non singular at $(r,t)=(1,0)$
and the curvature of the curve $G(\exp(t\sqrt{-1}))$ at $t=0$
is larger than $1/G(1)$.
In general, we can not judge only by the section at $r=1$. Indeed
$x^3+3x^2-x+1$ fulfills the condition but $-x^3+3x^2+x+2$ does not. 
They require a detailed study around $(r,t)=(1,\pi)$.
\bigskip

Irreducibility of lacunary polynomials is 
a classical subject and many related works are found in literature,
see for e.g. \cite{Ljunggren:60, Mills:85, Schinzel:00}. To make
explicit the constants $n_0$ in Theorem \ref{Gen}, 
the reader may consult \cite{Schinzel:67,Filaseta-Ford-Konyagin}. 

\bigskip

The set of simple Parry numbers is dense in $[1,\infty)$. We know little on 
the topology of the set of non-Parry Perron numbers in $\R$, neither on 
the set of their conjugates in $\C$.
\bigskip

\begin{center}
{\bf Acknowledgments}
\end{center}

This work is initiated by a question posed by Anne Bertrand, whether $\beta_4$
is a Perron number. The author is 
grateful to 
Shin'ichi Yasutomi and Kan Kaneko, for stimulating discussion, in particular, 
on the asymptotic behavior of roots of a parametrized family of polynomials. 

\bibliographystyle{amsplain}
\bibliography{../reflist}

\end{document}